\ifx\documentclass\undefined
\documentstyle[12pt]{article}
\else
\documentclass[12pt]{article}
\fi

\usepackage{color}

\usepackage{ascmac}

\newcommand{\beq}{\begin{eqnarray*}}
\newcommand{\eeq}{\end{eqnarray*}}
\makeatletter
\renewcommand{\theequation}{\thesection.\arabic{equation}}
\@addtoreset{equation}{section}
\def\eqnarray{%
\stepcounter{equation}%
\let\@currentlabel=\theequation
\global\@eqnswtrue
\global\@eqcnt\z@
\tabskip\@centering
\let\\=\@eqncr
$$\halign to \displaywidth\bgroup\@eqnsel\hskip\@centering
$\displaystyle\tabskip\z@{##}$&\global\@eqcnt\@ne
\hfil$\displaystyle{{}##{}}$\hfil
&\global\@eqcnt\tw@$\displaystyle\tabskip\z@{##}$\hfil
\tabskip\@centering&\llap{##}\tabskip\z@\cr}
\makeatother
\newtheorem{theorem}{Theorem}[section]
\newtheorem{lemma}[theorem]{Lemma}

\newsavebox{\toy}
\savebox{\toy}{\framebox[0.65em]{\rule{0cm}{1ex}}}
\newcommand{\QED}{\usebox{\toy}}
\def\nlni{\par\ifvmode\removelastskip\fi\vskip\baselineskip\noindent}

\begin{document}
%%%%%%% DOUBLE SPACED %%%%%%%%
\setlength{\baselineskip}{15pt}
\title{
A remark on conditions that a diffusion in the natural scale is a martingale
}
\author{
Yuuki Shimizu and Fumihiko Nakano
\thanks{
Department of Mathematics,
Gakushuin University,
1-5-1, Mejiro, Toshima-ku, Tokyo, 171-8588, Japan.
%
%e-mail : 
%fumihiko@math.gakushuin.ac.jp
}
}
%\date{最終更新日：}
\maketitle
%第一ページの番号を消す
%\thispagestyle{empty}
%%%%%%% ABSTRACT %%%%%%%%%%%%%
\begin{abstract}
We 
consider a diffusion processes 
$\{ X_t \}$ 
on an interval in the natural scale. 
Some results 
are known under which   
$\{ X_t \}$ 
is a martingale, and we give simple and analytic proofs for them.
% on the condition that a diffusion process in a natural scale is a martingale.
%
\end{abstract}


%\tableofcontents
%%%%% INTRODUCTION %%%%%%%%%%%%%%%%%%%%%%%%%%%%%%%%%
\section{Introduction}
Let 
$-\infty \le l_- < l_+ \le \infty$
and let 
$m$ 
be a Borel measure with  
$\mbox{ supp } m = (l_-, l_+)$.
We denote by 
$\Bigl\{
\{ X_t \}_{t \ge 0}, 
\{ P_x \}_{x \in (l_-, l_+) }
\Bigr\}$
the minimal diffusion process on 
$(l_-, l_+)$
with the speed measure 
$m$
and 
the scale function 
$s(x)= x$. 
It is well known that 
a local martingale 
$\{ X_t \}$ 
is a martingale if and only if 
$\{X_T : T\mbox{ is a stopping time with }\; T \le t\}$
is uniformly integrable for any 
$t \ge 0$. 
Here 
our aim is to have more explicit condition for the one-dimensional diffusions in the natural scale. 
If 
$| l_{\pm} | < \infty$, 
$\{ X_t \}$ 
is bounded so that it is a martingale.
If 
$l_- = - \infty$, 
$l_+ < \infty$, 
this can be reduced to the case of 
$l_- < \infty$, $l_+ = \infty$
by replacing 
$X_t$ 
by 
$-X_t$. 
Hence 
it suffices to consider the following two cases. 
\\

Case I : 
$-\infty< l_-$, $l_+ = +\infty$, 
\quad
Case II : 
$l_- = - \infty$, 
$l_+ = +\infty$.
\\
Let 
$P(l_-, l_+)$
be the set of Borel measures on 
$(l_-, l_+)$, 
and for 
$\mu \in P(l_-, l_+)$  
let 
$P_{\mu}(\cdot)
:=
\int_{(l_-,l_+)} P_x(\cdot) \mu(dx)$.
According to 
Lemma \ref{martin}
(\cite{kotani}, Lemma 2), 
$\{ X_t^{\tau} \}$
is a 
$P_{\mu}$-martingale 
for some 
$\mu \in P(l_-, l_+)$ 
with 
$\int_{(l_-, l_+)} |x| \mu(dx)< \infty$ 
if and only if 
$\{ X_t^{\tau} \}$
is 
$P_x$-martingale 
for any 
$x \in (l_-, l_+)$. 
We further set 
\beq 
\tau_a 
&:=&
\inf
\left\{ 
t \ge 0 \middle| X_t = a 
\right\}, 
\quad
\tau_{\pm} 
:=
\lim_{a \to l_{\pm}} \tau_a, 
\quad
\tau := \tau_+ \wedge \tau_-
\\
X^{\tau}_t 
&:=&
X_{t \wedge \tau}.
\eeq
Kotani \cite{kotani} 
showed the following theorem.
\begin{theorem}\label{kotani}
{\bf \cite{kotani}}
$\{X^{\tau}_{t} \}$
is a 
$P_x$-martingale 
for any 
$x\in (l_-, l_+)$
if and only if\\
Case I : 
\[ 
\int_{[r,l_+)} x m(dx) =\infty , 
\quad
r \in (l_-, \infty) 
\]
Case II : 
\[ 
\int_{[r,l_+)} x m(dx) =\infty
\;\mbox{ and }\;
\int_{(l_-,r]} |x| m(dx) = \infty,  
\quad
r \in (-\infty,\infty).
\]
\end{theorem}
By 
Feller's criterion, 
$P_x (\tau_{\sharp} = \infty) = 1$
if 
$| l_{\sharp} | = \infty$, 
$\sharp = \pm \infty$. 
Thus 
Theorem \ref{kotani} 
implies that 
$\{ X^{\tau}_t \}$
is a martingale if and only if the boundaries at infinity are natural.
Hulley, Platen 
\cite{hardy}
derived another condition.
Let 
\[
{\cal L}f 
:=
\frac {d^2}{d m dx}f
\]
be the generator of 
$\{ X_t \}$ 
and for 
$\lambda > 0$ 
let 
$f_-$
(resp. $f_+$)
be the positive increasing 
(resp. positive decreasing)
solution to the equation 
${\cal L} f = \lambda f$, 
which are unique up to constants unless the boundary is regular.
\begin{theorem}\label{hardy1}
{\bf \cite{hardy}}
$\{X^{\tau}_{t} \}$
is a 
$P_x$-martingale 
for any 
$x\in (l_-, l_+)$
if and only if\\
Case I : 
\[  
\lim_{z \to \infty} f'_-(z) = \infty  
\]
Case II : 
\[ 
\lim_{z \to \infty } f'_-(z) = \infty 
\;\mbox{ and }\;
\lim_{z \to -\infty} f'_+(z)= -\infty.
\]
\end{theorem} 
Gushchin, Urusov, and Zervos 
\cite{urusov}
derived a condition that 
$\{ X^{\tau}_t \}$ 
is a submartingale or a supermartingale.
\begin{theorem}
\label{submarurusov}
{\bf \cite{urusov}}
$\{ X^{\tau}_t \}$
is a 
$P_x$-submartingale 
if and only if 
$\int_r^{\infty} x m(dx) = \infty$, 
$r \in (l_-, l_+)$. 
\end{theorem}
By 
\cite{hardy} Proposition 3.16, 3.17, 
this condition is equivalent to 
$\lim_{t \to \infty} f'_- (t) = \infty$.
Together 
with Theorem \ref{submarurusov}
we thus have  
\begin{theorem}\label{submar1}
$\{ X^{\tau}_t \}$
is a 
$P_x$-submartingale 
if and only if 
$\lim_{t \to \infty}
f'_-(t) = \infty$.
\end{theorem}
Moreover in
\cite{urusov}, 
they further derived a condition in Case I such that 
$\{ X^{\tau}_{t} \}$
is a strict 
$P_x$ supermartingale, that is, 
$\{ X^{\tau}_{t} \}$
is a 
$P_x$-supermartingale but is not a 
$P_x$-martingale.
\begin{theorem}
\label{expectation}
{\bf \cite{urusov}}
Let 
$- \infty < l_-$, $l_+ = \infty$.
Then 
$\{ X^{\tau_-}_{t} \}$
is a strict 
$P_x$-supermartingale 
if and only if
\beq 
\lim_{t \to \infty} 
E_x [ X_{t \wedge \tau_-} ] = l_- 
\eeq
for any 
$x \in (l_-, l_+)$. 
\end{theorem}
We believe that 
Theorem \ref{expectation} 
is also true for 
$l_- = - \infty$.
The goal 
of this paper is : \\
(1)
To give 
a simple analytic proof of 
Theorem \ref{submar1} 
without using the results in 
\cite{hardy}.
We note that 
the proofs of 
Proposition 3.16, 3.17 in \cite{hardy} 
is more or less probabilistic using  Tanaka's formula.\\
(2)
To give 
a simple analytic proof of 
Theorem \ref{expectation} ; 
the original proof 
of that in 
\cite{urusov} 
is done by embedding 
$\{ X_t \}$ 
into the geometric Brownian motion on the torus.

The rest 
of this paper is organized as follows. 
In Section 2(resp. Section 3), 
we give a proof of Theorem \ref{submar1}
(resp. Theorem \ref{expectation}).
In Appendix, 
we prepare some tools for these proofs.
%%
%------ 証明１-----------
\section{A proof of Theorem 
\ref{submar1}}
In Case I, 
the statement follows from Theorem \ref{hardy1}, 
for 
$\{ X^{\tau_-}_t \}$
is always a 
$P_x$-supermartingale being bounded from below.
Henceforth 
we consider 
Case II. 
Suppose 
$\{ X_t \}$
is a $P_x$-submartingale and let 
$z < x$.
Then 
$\{ X^{\tau_z}_t \}$ 
is bounded from blow so that it is a 
$P_x$-martingale.
For
$\lambda > 0$, 
let 
$f^z_-$
(resp. $f^z_+$)
be the positive increasing 
(resp. positive decreasing) 
solution to the equation 
${\cal L} f = \lambda f$ 
such that 
$f^z_- (z) = 0$.
Then we have 
\[ 
f_-^z(x) = 
f_-(x) 
- 
\frac{f_-(z)}{f_+(z)} f_+(x),
\quad 
f_+^z(x) = f_+(x).
\]
Since 
$f_+'$
is increasing, we have
\beq
f'_-(x) &=& 
{f_-^{z}}'(x) + \frac{f_-(z)}{f_+(z)} f'_+(x) 
\ge 
{f_-^z}'(x) + \frac{f_-(z)}{f_+(z)} f'_+(z),
\quad
x \in (z, \infty).
\eeq
Applying 
Theorem \ref{hardy1} 
to 
$\{ X^{\tau_z}_t \}$ 
yields  
$\lim_{t \to \infty}{f_-^z}'(t)= \infty$
and thus 
$\lim_{t \to \infty} f'_-(t)=\infty$.
Conversely, suppose
$\lim_{t \to \infty}
f'_-(t) = \infty$ 
and let 
$z < x$. 
Then 
\beq
\lim_{z \to \infty} 
z 
\int^{\infty}_0 e^{-\lambda t} P_x( \tau_z < t ) dt 
&=&
\lim_{z \to \infty} 
\frac {z}{\lambda}
E_x[ e^{-\lambda \tau_z} ] 
= 
\lim_{z \to \infty} \frac{z}{\lambda} \frac{f_-(x)}{f_-(z)} 
\\
&=& 
\lim_{z \to \infty} 
\frac{f_-(x)}{\lambda} 
\frac{1}{f'_-(z)} 
= 0 
\eeq
where 
we used Lemma \ref{hitting} 
and l'Hospital's rule.
By Fatou's lemma,
\[ 
\int^{\infty}_0 e^{-\lambda t} \liminf_{z \to \infty} zP_x( \tau_z < t ) dt = 0. 
\]
Hence 
$\liminf_{z \to \infty} 
zP_x(\tau_z < t) = 0$ 
so that we can find a sequence 
$\{z_n\} \subset (x, \infty)$
with  
$\lim_{n \to \infty}z_n = \infty$
such that 
\[ 
\lim_{n \to \infty} z_n P_x ( \tau_{z_n} < t ) = 0.
\]
On the other hand 
$\{ X^{\tau_{z_n}}_{t} \}$
is a 
$P_x$-submartingale being bounded from above and 
\[ 
x \le E_x [ X_{t \wedge \tau_{z_n}} ]  
= 
z_n
P_x(\tau_{z_n} < t ) + 
E_x [ X_t ; \tau_{z_n} \ge t ]. 
\]
Since 
$\lim_{n \to \infty}
P_x(\tau_{z_n} \ge t ) =1$, 
$x \le E_x [X_t]$. 
Markov property 
implies  
$\{ X_t \}$
is a $P_x$-submartingale.
\QED
%
%
%-----定理漸近挙動の証明--------
\section{A proof of Theorem \ref{expectation}}
Without losing generality, 
we may suppose 
$l_- < 0$.
For
$\lambda > 0$, 
let 
$f_-$
(resp. $f_+$)
be the positive increasing 
(resp. positive decreasing) 
solution to the equation 
${\cal L} f = \lambda f$ 
such that 
$f_- (l_-) = 0$.
Let 
$G$ 
be Green's function of 
${\cal L}$ : 
\beq
G(x,y,\lambda)
&:=&
\left\{
\begin{array}{cc}
\frac{1}{h} f_-(y) f_+(x)
& 
(y < x ) \\
\frac{1}{h}f_-(x) f_+(y)
&
(x \le y)
\end{array}
\right. 
\\
h
&:=&
f_+(x) f'_-(x)
-
f_-(x) f'_+(x).
\eeq
Then we have
\begin{equation}
\int^{\infty}_{l_-} 
G(x,y,\lambda) (y-l_-) m(dy) 
= 
E_x 
\Big[ 
\int^{\infty}_0 
e^{-\lambda t} 
(X_{t \wedge \tau_-} -l_-) dt 
\Big].
\label{Green}
\end{equation}
Let 
$\alpha_+ := 
\lim_{t \to \infty} f_+(t)$.
Then  
$f_+' \in L^1(a, \infty)$ 
for 
$a \in (l_-, \infty)$ 
and 
\[ 
f_+(x)= 
\alpha_+ - 
\int^{\infty}_x f'_+ (y) dy.
\]
Therefore
$\lim_{x \to \infty } 
f'_+ (x)=0$. 
The equation 
${\cal L}f_+ = \lambda f_+$ 
yields
\beq
f'_+(x) 
&=&
-
\lambda 
\int^{\infty}_x 
f_+(y) m(dy)
\\
f_+(x) 
&=& 
\alpha_+ 
+ 
\lambda 
\int^{\infty}_x
(y-x)f_+(y)m(dy)
\eeq
so that we have
\beq
\lambda \int^{\infty}_{x} y f_+(y) m(dy)
&=&
f_+(x) - \alpha_+ -x f'_+(x).
\eeq
Similarly, 
\beq
f'_- (y)
&=&
f'_- (l_-)
+
\lambda
\int_{l_-}^y f_- (z) m(dz)
\\
f_- (x)
&=&
f_-' (l_-) (x - l_-)
+
\lambda
\int_{l_-}^x 
(x - y) f_- (y) m(dy)
\\
\lambda
\int_{l_-}^x 
y f_- (y) m(dy)
&=&
f_-' (l_-) (x - l_-)  - f_- (x)
+
\lambda
x 
\int_{l_-}^x
f_- (y) m(dy).
\eeq
Substituting them into 
(\ref{Green}) 
yields
\begin{eqnarray}
\int^{\infty}_0 
e^{-\lambda t} 
E_x[ X_{t \wedge \tau_-} -l_-] dt 
&=&
\frac{x - l_-}{\lambda}
-
\frac{\alpha_+ f_-(x )}
{\lambda h}.
\label{star}
\end{eqnarray}
We note that 
(\ref{star})
and 
Lemma \ref{martin}
also proves Theorem \ref{kotani} in 
Case I.\\

\noindent
Suppose 
$\{X^{\tau_-}_{t} \}$ 
is a strict $P_x$-supermartingale.
The discussion above implies 
$\alpha_+ > 0$.
We shall show below that  
\begin{equation}
\lim_{\lambda \to 0}  
\Big( 
x - l_-
- 
\frac{\alpha_+ f_-(x)}{h} \Big) = 0. 
\label{lambdalimit}
\end{equation}
Let 
$\phi$, $\psi$ 
be the solution to 
${\cal L} f = \lambda f$ 
with the initial condition
\beq
&&
\phi(0) = 1, \quad \phi'(0) = 0
\\
&&
\psi(0) = 0, \quad \psi'(0) = 1.
\eeq
Then 
$f_{\pm}$
satisfy
\beq
&&
f_+(x) = \phi(x) - \left( \lim_{x \to \infty}  \frac{\phi(x)}{\psi(x)} \right) \psi(x),
\quad 
f_-(x) = \phi(x) - \left( \lim_{x \to l_-} \frac{\phi(x) }{\psi(x) } \right) \psi(x).
\eeq
$\psi$, $\psi$ 
can be composed by the method of successive approximation : 
\beq
&&
\phi(x) = 1 + \sum_{n = 1} ^{\infty} \lambda^n \phi_n(x),
\quad 
\phi_0(x) = 1,
\quad 
\phi_n(x) = \int^x_0 (x-y) \phi_{n-1}(y) m(dy) 
\\
&&
\psi(x) = x + \sum_{n = 1}^{\infty} \lambda^n \psi_n(x),
\quad 
\psi_0(x) = x,
\quad 
\psi_n(x) = \int^{x}_{0} (x-y) \psi_{n-1} (x) m(dy)
\eeq
which is convergent locally uniformly w.r.t. 
$\lambda$
\cite{ito}
which yields
\[ 
\lim_{\lambda \to 0} \phi(x) = 1,
\quad 
\lim_{\lambda \to 0} \phi'(x) = 0, 
\quad 
\lim_{\lambda \to 0} \psi(x) = x,
\quad 
\lim_{\lambda \to 0}\psi'(x) = 1.
\]
Moreover
\[ 
\lim_{\lambda \to 0} 
\left( -\lim_{x \to \l_-} \frac{\psi(x)}{\phi(x)}
\right)  
= 
\lim_{\lambda \to 0} \left(  
\int^{0}_{l_-} 
\frac{1}{(\phi(x))^2}dx 
\right) =  \int^0_{l_-} dx = - l_- 
\]
implies
\[ 
\lim_{\lambda \to 0} 
f_-(x) 
= 
1 - \frac{x} {l_-}, 
\quad 
\lim_{\lambda \to 0} f'_-(x) = -\frac{1}{l_-}.
\]
On the other hand, 
by 
$\alpha_+ > 0$ 
and by 
Lemma \ref{point1}, 
we have 
$\int^{\infty}_r xm(dx) < \infty$, 
$r \in (l_-,\infty)$
so that we can find
$g$ 
satisfying 
\[ 
g(x) = 1 + \lambda \int^{\infty}_x (y-x) g(y)  m(dy) 
\]
by successive approximation.
Using 
$\alpha_+ > 0$, 
$\lim_{t \to \infty} f'_+(t)=0$, 
$\lim_{t \to \infty} g(t)=1$ 
and 
$\lim_{t \to \infty} g'(t) = 0$, 
we have 
\[ 
f_+(x) g'(x) - f_+'(x) g(x) = 0
\]
which implies  
$f_+(x) = 
C g(x)$
for some positive constant
$C$. 
Because  
$\lim_{\lambda \to 0}
g(x) = 1$, 
$\lim_{\lambda \to 0}
g'(x) = 0$, 
\[ 
\lim_{\lambda \to 0} f_+(x) = 
C, 
\quad
\lim_{\lambda \to 0}
f'_+(x) = 0. 
\]
Therefore
\[ 
\lim_{\lambda \to 0}  
\left( 
x - l_-
- 
\frac{\alpha_+ f_-(x)}{h} 
\right) 
= 
x-l_- - 
\frac{
C
\left( 1 - \frac{x}{l_-}\right)
}
{
C \cdot
\left(
\frac{-1}{l_-}
\right) 
- 
0 \cdot 
\left( 1 - \frac{x}{l_-}\right) 
}  
= 0
\]
proving
(\ref{lambdalimit}).
Since 
$X_{t \wedge \tau_-}$ 
is a supermartingale, 
$f(t) := E_x [ X_{t \wedge \tau_-} - l_-] \in C^1[0, \infty)$
is monotone decreasing
which shows that 
$\lim_{t \to \infty}f(t)$ 
exists and 
$f' \in L^1(0, \infty)$. 
Thus by 
(\ref{star}) and 
Lemma \ref{Laplace}
\[ 
\lim_{t \to \infty} 
E_x [ X_{t \wedge \tau_-} - l_-] 
= 0.
\]
Conversely, suppose that 
$\displaystyle 
\lim_{t \to \infty} 
E_x[X_{t \wedge \tau_-} - l_- ] = 0$.
Then 
\[ 
\lim_{\lambda \to 0} 
\lambda 
\int^{\infty}_0 
e^{-\lambda t} 
E_x[X_{t \wedge \tau_- } - l_-] dt 
= 0 
\]
which implies 
$\alpha_+ > 0$ 
since otherwise it would contradict to 
(\ref{star}), 
(\ref{lambdalimit}).
Therefore 
$\{ X^{\tau_-}_t \}$ 
is not a martingale.
\QED
%
%--付録-----------------
\section{Appendix}
\begin{lemma}\label{martin}
{\bf (Lemma 2 in \cite{kotani})}\\
Suppose 
$\{X_{t \wedge \tau_-} \}$
is a $P_{\mu}$-martingale for some 
$\mu \in P\big(l_-,\infty \big)$.
Then 
for any 
$t \ge 0$, $x \in (l_-, \infty)$, 
\begin{equation}
E_x [ X_{t \wedge \tau_-} ]
=
x.
\label{x}
\end{equation}
Conversely, if 
(\ref{x}) 
is valid, then 
$\{X_{t \wedge \tau_-} \}$
is a $P_{\mu}$-martingale for any  
$\mu \in P\big(l_-,\infty \big)$ 
with 
$\int_{l_-}^{\infty} |x| \mu (dx) < \infty$. 
\end{lemma}
\begin{lemma}
\label{point1}
Let 
$\lambda >0$
and let 
$f_+$ 
be the positive decreasing solution to 
${\cal L}f = \lambda f$ 
with 
$\displaystyle \alpha _+ :=\lim_{x \to \infty} f_+(x)$.
Then 
the following three conditions are equivalent. 
\beq
(1)&&\quad
\alpha_+=0 
\\
(2)&&\quad
%\Longleftrightarrow 
\int^{\infty}_a y m(dy)= \infty 
\\
(3)&& \quad
%\Longleftrightarrow  
\lambda \int ^{\infty}_x  (y-x) f_+(y) m(dy) = f_+(x).
\eeq
\end{lemma}
\begin{lemma}\label{hitting}
Let 
$f_{\pm}$ 
be the ones defined in the proof of 
Theorem \ref{expectation}. 
Then 
\beq
&&
E_x [ e^{-\lambda \tau_a } ]
= 
\frac{f_+(x)}{f_+(a)},
\quad
a < x 
\\
&&
E_x [e^{-\lambda \tau_b} : \tau_b < \tau_-] = \frac{f_-(x)}{f_-(b)},
\quad
- \infty \le l_- < x < b. 
\eeq
\end{lemma}
\begin{lemma}
\label{Laplace}
Suppose 
$f \in C^1[0, \infty)$
and 
$f' \in L^1(0, \infty)$.
Then 
\\
(1)
$\lim_{t \to \infty} f(t)$
exists, and 
\\
(2)
$
\lim_{t \to \infty} f(t)
=
\lambda
\lim_{\lambda \downarrow 0}
\int_0^{\infty}
e^{- \lambda t} f(t) dt.
$
\end{lemma}
%

%\vspace*{1em}
%\noindent {\bf Acknowledgement }
%This work is partially supported by 
%Grant-in-Aid for Scientific Research (C) no.26400145.

%%%%% REFERENCES %%%%%%%%%%%%%%%%%%%%%
%

%
\end{document}